# NONLINEAR SCHRÖDINGER EQUATIONS WITH REPULSIVE HARMONIC POTENTIAL AND APPLICATIONS

RÉMI CARLES

ABSTRACT. We study the Cauchy problem for Schrödinger equations with repulsive quadratic potential and power-like nonlinearity. The local problem is well-posed in the same space as that used when a confining harmonic potential is involved. For a defocusing nonlinearity, it is globally well-posed, and a scattering theory is available, with no long range effect for any superlinear nonlinearity. When the nonlinearity is focusing, we prove that choosing the harmonic potential sufficiently strong prevents blow-up in finite time. Thanks to quadratic potentials, we provide a method to anticipate, delay, or prevent wave collapse; this mechanism is explicit for critical nonlinearity.

1. INTRODUCTION

Consider the Schrödinger equation

$$(1.1) \qquad i\partial_t u + \frac{1}{2}\Delta u = V(x)u + \lambda|u|^{2\sigma}u, \quad (t,x) \in \mathbb{R} \times \mathbb{R}^n,$$

with $\sigma > 0$, $\sigma < 2/(n-2)$ if $n \geq 3$, $\lambda \in \mathbb{R}$ and $V$ is a real-valued potential $V : \mathbb{R}^n \to \mathbb{R}$. If $V \in L^\infty + L^p$, for some $p \geq 1$, $p > n/2$, then the Cauchy problem in $H^1(\mathbb{R}^n)$ associated to (1.1) is known to be locally well-posed; it may also be globally well-posed or lead to blow-up in finite time (see e.g. [Caz93]).

If the potential is smooth, $V \in C^\infty$, non-negative, and if its derivatives of order at least two are bounded, then the same holds in the domain of $\sqrt{-\Delta + V}$ (see [Oh89], [Caz93]). When $n = 1$ and $V$ is non-negative with superquadratic growth, then the fundamental solution for (1.1) with $\lambda = 0$ is nowhere $C^1$ ([Yaj96]), but smoothing properties make it possible to solve the nonlinear problem (1.1) in some cases ([YZ01]).

If $V$ is non-positive, then $-\Delta + V$ is essentially self-adjoint on $C_0^\infty(\mathbb{R}^n)$ provided that there exist some constants $a, b$ such that $V(x) \geq -a|x|^2 - b$ (see [RS75], p. 199). If $-V$ has superquadratic growth, then it is not possible to define $e^{-it(-\Delta+V)}$ (see [DS63], Chapter 13, Sect. 6, Cor. 22). In this paper, we study the Cauchy problem

$$(1.2) \qquad \begin{cases} i\partial_t u + \frac{1}{2}\Delta u = -\omega^2 \frac{|x|^2}{2}u + \lambda|u|^{2\sigma}u, \quad (t,x) \in \mathbb{R} \times \mathbb{R}^n, \\ u_{|t=0} = u_0, \end{cases}$$

with $\omega, \sigma > 0$, $\sigma < 2/(n-2)$ if $n \geq 3$, $\lambda \in \mathbb{R}$, and

$$u_0 \in \Sigma := \left\{ f \in H^1(\mathbb{R}^n) \, ; \, |x|f \in L^2(\mathbb{R}^n) \right\}.$$

The Hilbert space $\Sigma$ is equipped with the norm

$$\|f\|_\Sigma = \|f\|_{L^2} + \|\nabla f\|_{L^2} + \|xf\|_{L^2}.$$

Another motivation for studying (1.2) lies in the study of finite time blow up for the Cauchy problem

$$
(1.3) \quad \begin{cases} i\partial_t u + \frac{1}{2}\Delta u = \lambda |u|^{2\sigma} u, & (t,x) \in \mathbb{R} \times \mathbb{R}^n, \\ u_{|t=0} = u_0. \end{cases}
$$

It is well-known that if $u_0 \in \Sigma$, $\lambda < 0$, $\sigma \geq 2/n$ and

$$
(1.4) \quad E(u_0) := \frac{1}{2}\|\nabla u_0\|_{L^2}^2 + \frac{\lambda}{\sigma+1}\|u_0\|_{L^{2\sigma+2}}^{2\sigma+2} < 0,
$$

then $u$ blows up in finite time, that is there exists $T > 0$ such that

$$
\lim_{t \to T} \|\nabla_x u(t)\|_{L^2} = +\infty.
$$

This is proven by the general approach of Zakharov-Glassey ([Gla77], [Caz93]). There exist numerical evidences suggesting that the introduction of a stochastic white noise in (1.3) may amplify or prevent blow-up formation (see [DDM02]); in our paper, we confine ourselves to a deterministic framework.

It is shown in [CW92] that if the initial datum $u_0$ is replaced by $u_0(x)e^{-ib|x|^2}$ with $b > 0$ sufficiently large, then the blow-up time is anticipated (and is $O(b^{-1})$). On the other hand, if $u_0$ is replaced by $u_0(x)e^{ib|x|^2}$ with $b > 0$ sufficiently large, then no blow-up occurs.

Our approach is suggested by the semi-classical régime for linear Schrödinger equation with potential. Consider the initial value problem

$$
\begin{cases} i\varepsilon \partial_t u^\varepsilon + \frac{1}{2}\varepsilon^2 \Delta u^\varepsilon = V(x) u^\varepsilon, & (t,x) \in \mathbb{R} \times \mathbb{R}^n, \\ u^\varepsilon_{|t=0} = u_0^\varepsilon, \end{cases}
$$

where $V \in C^\infty(\mathbb{R}^n, \mathbb{R})$, $\varepsilon \in {]0,1]}$. In the semi-classical limit $\varepsilon \to 0$, the energy of the solution $u^\varepsilon$ is carried by bicharacteristics, which are the integral curves associated to the classical Hamiltonian

$$
p(t,x,\tau,\xi) = \tau + \frac{|\xi|^2}{2} + V(x).
$$

If the energy tends to concentrate in this case, one can expect that for (1.1) with a focusing nonlinearity ($\lambda < 0$), blow-up in finite time, which corresponds to the concentration of the mass, may occur. Bicharacteristic curves solve

$$
\begin{cases} \dot t = 1, \\ \dot x = \xi, \\ \dot \tau = 0, \\ \dot \xi = -\nabla V(x). \end{cases}
$$

Rays of geometric optics, which are the projection of bicharacteristic curves on $(t,x)$ space, are of the form $x = x(t)$, with

$$
(1.5) \quad \begin{cases} \ddot x + \nabla V(x) = 0, \\ x(0) = x_0, \ \dot x(0) = \xi_0. \end{cases}
$$

If the initial datum is of the form $u_0^\varepsilon(x) = f(x)e^{i\varphi(x)/\varepsilon}$, then $\xi_0 = \nabla \varphi(x_0)$. We give three examples which are at the origin of this work, and which correspond to cases where (1.5) can easily be solved.

*Example* 1. Suppose $V \equiv 0$. Then the solutions of (1.5) are

$$
x(t) = x_0 + t\nabla \varphi(x_0).
$$

If no oscillation is present in the initial datum, rays are parallel. More interesting is the case of quadratic oscillations (see also [CFG]). If $\varphi(x) = -b|x|^2/2$ with $b > 0$, then rays are given by $x(t) = x_0(1 - bt)$, and meet at the origin at time $t = 1/b$ (see Figure 1). There is focusing, which suggests that in a nonlinear setting, such oscillations may cause wave collapse. If $\varphi(x) = b|x|^2/2$ with $b > 0$, then rays are given by $x(t) = x_0(1 + bt)$, and met at the origin at time $t = -1/b$ (in the past). In particular, they are spread out for positive times, which suggests that in a nonlinear setting, such oscillations may prevent wave collapse. The intuition described on the last two cases is confirmed by the results of Cazenave and Weissler [CW92].

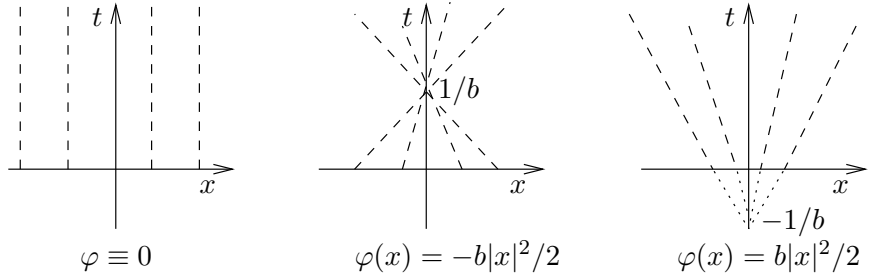

FIGURE 1. Geometry of rays: case $V \equiv 0$.

*Example* 2. Suppose $V(x) = \omega^2 |x|^2/2$, with $\omega > 0$. In the case $\varphi \equiv 0$, rays are given by $x(t) = x_0 \cos(\omega t)$, and meet at the origin at time $t = \pi/(2\omega)$ (see Figure 2). The first example suggests that blow-up may happen more easily than when $V \equiv 0$. This phenomenon is reinforced by the case $\varphi(x) = -\omega \tan(\omega t_0)|x|^2/2$, $|t_0| < \pi/(2\omega)$, where
$$x(t) = \frac{x_0}{\cos(\omega t_0)} \cos \omega(t + t_0).$$
Rays meet at the origin at time $t = \pi/(2\omega) - t_0$. If $t_0 > 0$, focusing is anticipated, while if $t_0 < 0$, it is delayed (but in no case prevented). This geometry is to be compared with the second case of the first example.

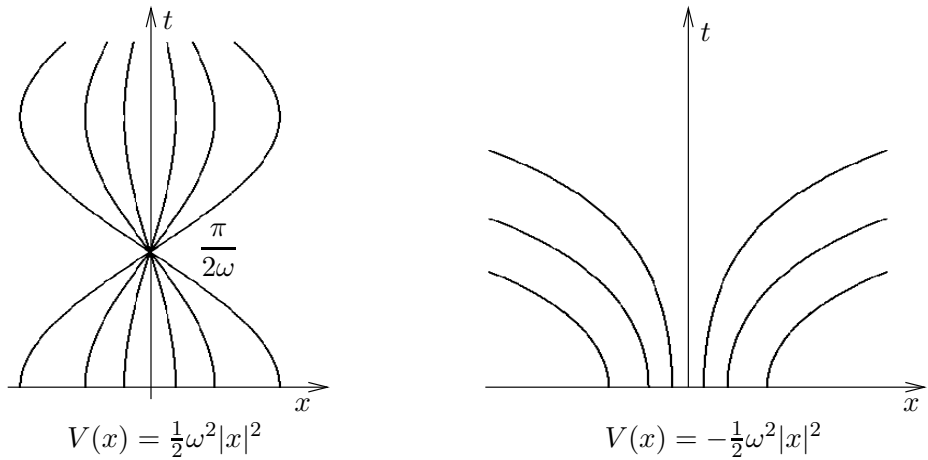

FIGURE 2. Geometry of rays, with $\varphi \equiv 0$.

*Example* 3. Suppose $V(x) = -\omega^2 |x|^2/2$, with $\omega > 0$. In the case $\varphi \equiv 0$, rays are given by $x(t) = x_0 \cosh(\omega t)$, and are strongly dispersed for positive times (see Figure 2). This geometry is to be compared with the third case of the first example; rays are scattered, but go to infinity exponentially fast, instead of algebraically. If $\varphi(x) = -\omega b |x|^2/2$ with $b > 0$, then rays are given by $x(t) = x_0(\cosh(\omega t) - b\sinh(\omega t))$, and their behavior is far less singular than in the case $V \equiv 0$. This is a first hint that such potentials may prevent blow-up.

Inspired by the second example, we proved in [Car02b] that if $u$ solves

$$(1.6) \quad \begin{cases} i\partial_t u + \frac{1}{2}\Delta u = \omega^2 \frac{|x|^2}{2} u + \lambda |u|^{2\sigma} u, & (t,x) \in \mathbb{R} \times \mathbb{R}^n, \\ u_{|t=0} = u_0, \end{cases}$$

then under Condition (1.4) (as a matter of fact, the value $E(u_0) = 0$ is allowed), $u$ blows up at time $T \leq \pi/(2\omega)$; choosing $\omega$ large enough, the blow-up time is therefore anticipated by the action of a confining magnetic field. This is intimately related to the dynamics of the linear Schrödinger equation with a confining harmonic potential, whose solution is given by Mehler's formula (see e.g. [FH65]), for $|t| < \pi/(2\omega)$,

$$(1.7) \quad u(t,x) = e^{-in\frac{\pi}{4}\operatorname{sgn} t} \left| \frac{\omega}{2\pi \sin \omega t} \right|^{n/2} \int_{\mathbb{R}^n} e^{\frac{i\omega}{\sin \omega t}\left(\frac{x^2+y^2}{2}\cos \omega t - x \cdot y\right)} u_0(y) dy.$$

At time $t = \pi/(2\omega)$, the fundamental solution is singular. In the nonlinear case, two phenomena cumulate: because of the linear dynamics, the solution tends to focus near the origin as time goes to $\pi/(2\omega)$; when the solution is sufficiently concentrated, the nonlinear term becomes important and causes wave collapse.

Setting $\lambda = 0$ in (1.2), we have the analog of Mehler's formula,

$$(1.8) \quad u(t,x) = e^{-in\frac{\pi}{4}\operatorname{sgn} t} \left| \frac{\omega}{2\pi \sinh \omega t} \right|^{n/2} \int_{\mathbb{R}^n} e^{\frac{i\omega}{\sinh \omega t}\left(\frac{x^2+y^2}{2}\cosh \omega t - x \cdot y\right)} u_0(y) dy.$$

Not only the kernel of $U_\omega(t) := \exp\{-it/2(-\Delta - \omega^2|x|^2)\}$ given by the above formula has no singularities for $t > 0$, but the dispersive effects are much stronger than in the case with no potential ($\omega = 0$). One might think that in the nonlinear case, the free dynamics can prevent the nonlinear mechanism of blow-up, at least for $\omega$ large. We prove that this holds true in Sect. 4. These results are summarized in the following theorem, whose first point is a consequence of [Car02b].

**Theorem 1.1.** *Let $u_0 \in \Sigma$, $\lambda < 0$, $\sigma \geq 2/n$ and $\sigma < 2/(n-2)$ if $n \geq 3$.*
*1. Assume that the solution $u$ to (1.3) has negative energy, that is (1.4) holds. Then $u$ blows up in finite time $T$. Let $\omega_0 = \pi/(2T)$; for any $\omega \geq \omega_0$, the solution to (1.6) blows up before time $\pi/(2\omega)$, and in particular, before time $T$.*
*2. If the initial datum $u_0$ satisfies*

$$\frac{1}{2}\|\nabla u_0\|_{L^2}^2 + \frac{\lambda}{\sigma+1}\|u_0\|_{L^{2\sigma+2}}^{2\sigma+2} < -\frac{\omega^2}{2}\|xu_0\|_{L^2}^2,$$

*then the solution to (1.2) blows up in finite time, in the future* or *in the past.*
*3. If the initial datum $u_0$ satisfies*

$$\frac{1}{2}\|\nabla u_0\|_{L^2}^2 + \frac{\lambda}{\sigma+1}\|u_0\|_{L^{2\sigma+2}}^{2\sigma+2} < -\frac{\omega^2}{2}\|xu_0\|_{L^2}^2 - \omega \left|\operatorname{Im} \int \overline{u_0} x \cdot \nabla_x u_0\right|,$$

*then the solution to (1.2) blows up in finite time, in the future* and *in the past.*
*4. There exists $\omega_1 > 0$ such that for any $\omega \geq \omega_1$, (1.2) has a unique, global, solution $u \in C(\mathbb{R}, \Sigma)$.*

In the particular case $\sigma = 2/n$, a change of variables relating the solutions of (1.3) to those of (1.6) or (1.2) shows explicitly how the blow-up for solutions to (1.3) can be anticipated, delayed or prevented by the introduction of quadratic potentials.

**Theorem 1.2.** *Let $u_0 \in \Sigma$, $\lambda < 0$, $\sigma = 2/n$.*
*Assume that the solution to (1.3) blows up at time $T > 0$.*

- *For any $\omega > 0$, the solution to (1.6) blows up at time $\arctan(\omega T)/\omega < T$.*
- *If $\omega < 1/T$, then the solution to (1.2) blows up at time $\arg\tanh(\omega T)/\omega > T$.*
- *If $\omega \geq 1/T$, then (1.2) has a unique, global, solution in $C(\mathbb{R}_+, \Sigma)$.*

*Remark.* As we recall in Section 2 (see Lemma 2.2), (1.8) provides Strichartz inequalities, that make it possible to study (1.2) with $u_0 \in L^2(\mathbb{R}^n)$ only, for $\sigma < 2/n$, or $u_0 \in L^2(\mathbb{R}^n)$ and $\|u_0\|_{L^2}$ small if $\sigma = 2/n$. Our goal is precisely to understand the other cases, that is why we shall always assume $u_0 \in \Sigma$.

*Remark.* Replacing $\omega$ with $\pm i\omega$ formally turns (1.6) into (1.2) (and *vice versa*), and (1.7) into (1.8). All the algebraic manipulations we perform in Section 2 can be retrieved using this argument; in particular, (2.3), Lemma 2.3 and the evolution law (2.10) can be deduced from the formulae established in [Car02b].

When $\lambda > 0$, solutions to (1.3) are known to be global, and the classical issue is to understand their asymptotic behavior as $t \to \pm\infty$. For $\sigma$ sufficiently large, the solutions are asymptotically free, while for $\sigma \leq 1/n$, a long range scattering theory is needed (see [Bar84], [Str74], [Str81], [Oza91], [HN98], [Car01]). Notice that for the initial value problem (1.2), it is not obvious that the formal conservations of mass and energy imply global existence in $\Sigma$ once local existence is known. These conservations read

$$\|u(t)\|_{L^2} \equiv \|u_0\|_{L^2} \ ; \ \frac{1}{2}\|\nabla_x u(t)\|_{L^2}^2 - \frac{\omega^2}{2}\|xu(t)\|_{L^2}^2 + \frac{\lambda}{\sigma+1}\|u(t)\|_{L^{2\sigma+2}}^{2\sigma+2} = cst.$$

For (1.6), the analog of these two conservation laws yields global existence in $\Sigma$ when $\lambda > 0$, for the energy is the sum of three positive terms. For (1.2), the energy functional is not always positive, even if $\lambda > 0$. We prove in Section 2 that a refined analysis of the conservation of energy, consisting in splitting the energy into two parts, yields global existence as soon as $\lambda > 0$ (and in other cases). Moreover, the strong dispersive properties of $U_\omega$ lead to a different scattering theory for (1.2); the nonlinearity $u \mapsto |u|^{2\sigma} u$ is always short range.

**Theorem 1.3.** *Let $\lambda, \sigma > 0$, with $\sigma < 2/(n-2)$ if $n \geq 3$.*
*1. For every $u_- \in \Sigma$, there exists a unique $u_0 \in \Sigma$ such that the maximal solution $u \in C(\mathbb{R}, \Sigma)$ to (1.2) satisfies*

$$\|U_\omega(-t)u(t) - u_-\|_\Sigma \underset{t \to -\infty}{\longrightarrow} 0.$$

*2. For every $u_0 \in \Sigma$, there exists a unique $u_+ \in \Sigma$ such that the maximal solution $u \in C(\mathbb{R}, \Sigma)$ to (1.2) satisfies*

$$\|U_\omega(-t)u(t) - u_+\|_\Sigma \underset{t \to +\infty}{\longrightarrow} 0.$$

This paper is organized as follows. In Section 2, we study the Cauchy problem for (1.2); we prove that it is locally well-posed in $\Sigma$, and in some cases, globally well-posed. This is the most important part of the paper. We analyze the corresponding scattering theory in Section 3, and prove Theorem 1.3. Theorem 1.1 is proven in Section 4, Theorem 1.2 in Section 5.

## 2. Solving the Cauchy problem

To solve the local Cauchy problem, we first introduce the classical notions of admissible pairs and Strichartz estimates.

**Definition 2.1.** A pair $(q,r)$ is **admissible** if $2 \leq r < \frac{2n}{n-2}$ (resp. $2 \leq r \leq \infty$ if $n=1$, $2 \leq r < \infty$ if $n=2$) and
$$\frac{2}{q} = \delta(r) := n\left(\frac{1}{2} - \frac{1}{r}\right).$$

Recall that $U_\omega(t)$ denotes the semi-group $\exp(it/2(\Delta + \omega^2|x|^2))$, which is given explicitly by (1.8).

**Lemma 2.2** (Strichartz estimates for $U_\omega$). *Let $\omega > 0$.*
*1. For any admissible pair $(q,r)$, there exists $C_r$ independent of $\omega > 0$ such that,*
$$\|U_\omega(.)\varphi\|_{L^q(\mathbb{R};L^r)} \leq C_r \|\varphi\|_{L^2}, \tag{2.1}$$
*for every $\varphi \in L^2(\mathbb{R}^n)$.*
*2. For any admissible pairs $(q_1, r_1)$ and $(q_2, r_2)$ and any interval $I$, there exists $C_{r_1,r_2}$ independent of $\omega > 0$ and $I$ such that*
$$\left\|\int_{I\cap\{s\leq t\}} U_\omega(t-s)F(s)ds\right\|_{L^{q_1}(I;L^{r_1})} \leq C_{r_1,r_2} \|F\|_{L^{q_2'}(I;L^{r_2'})}, \tag{2.2}$$
*for every $F \in L^{q_2'}(I; L^{r_2'})$.*

*Proof.* The semi-group $U_\omega$ is isometric on $L^2(\mathbb{R}^n)$, and from (1.8), it satisfies, for any $t \neq 0$ and $f \in L^1(\mathbb{R}^n)$,
$$\|U_\omega(t)f\|_{L^\infty} \leq \frac{1}{|2\pi t|^{n/2}} \|f\|_{L^1}.$$
It follows from the results of [KT98] that $U_\omega$ satisfies such Strichartz estimates as stated above. One can choose constants independent of $\omega > 0$ because the above dispersion estimate is independent of $\omega > 0$. One can actually take the same constants as in the case with no potential, $\omega = 0$. □

As mentioned in the introduction, this lemma makes it possible to study (1.2) if $u_0 \in L^2(\mathbb{R}^n)$ and $\sigma < 2/n$, by just mimicking the proof of the corresponding result for (1.3). Since our interest is to study (1.2) when $\sigma \geq 2/n$ to analyze finite time blow up, we assume $u_0 \in \Sigma$, and introduce two operators,
$$J(t) := \omega x \sinh(\omega t) + i\cosh(\omega t)\nabla_x \; ; \; H(t) := x\cosh(\omega t) + i\frac{\sinh(\omega t)}{\omega}\nabla_x. \tag{2.3}$$

For $f \in L^2(\mathbb{R}^n)$ and $t \in \mathbb{R}$, the property $J(t)f, H(t)f \in L^2(\mathbb{R}^n)$ implies $f \in \Sigma$:
$$i\nabla_x = \cosh(\omega t)J(t) - \omega\sinh(\omega t)H(t) \; ; \; x = \cosh(\omega t)H(t) - \frac{\sinh(\omega t)}{\omega}J(t). \tag{2.4}$$

These two operators are the formal analogs of those we used in [Car02b] to study (1.6), when $\omega$ is replaced by $\pm i\omega$. They have the remarkable property to be both Heisenberg observables *and* conjugate to $\nabla_x$ by a unitary factor.

**Lemma 2.3.** *The operators $J$ and $H$ satisfy the following properties.*
*1. They are Heisenberg observables,*
$$J(t) = U_\omega(t)i\nabla_x U_\omega(-t) \; ; \; H(t) = U_\omega(t)xU_\omega(-t), \tag{2.5}$$
*and consequently commute with the linear part of (1.2),*
$$\left[i\partial_t + \frac{1}{2}\Delta + \omega^2\frac{|x|^2}{2}, J(t)\right] = \left[i\partial_t + \frac{1}{2}\Delta + \omega^2\frac{|x|^2}{2}, H(t)\right] = 0.$$

2. They can be factorized as follows, for $t \neq 0$,

(2.6)
$$J(t) = i\cosh(\omega t)e^{i\omega \frac{|x|^2}{2}\tanh(\omega t)}\nabla_x \left(e^{-i\omega \frac{|x|^2}{2}\tanh(\omega t)} \cdot \right),$$
$$H(t) = i\frac{\sinh(\omega t)}{\omega}e^{i\omega \frac{|x|^2}{2}\coth(\omega t)}\nabla_x \left(e^{-i\omega \frac{|x|^2}{2}\coth(\omega t)} \cdot \right).$$

3. They yield modified Gagliardo-Nirenberg inequalities. Recall that if $r \geq 2$, with $r < 2n/(n-2)$ if $n \geq 3$, there exists $c_r$ such that for any $f \in H^1(\mathbb{R}^n)$,

$$\|f\|_{L^r} \leq c_r \|f\|_{L^2}^{1-\delta(r)}\|\nabla f\|_{L^2}^{\delta(r)}.$$

Then for every $f \in \Sigma$,

(2.7)
$$\|f\|_{L^r} \leq \frac{c_r}{(\cosh(\omega t))^{\delta(r)}}\|f\|_{L^2}^{1-\delta(r)}\|J(t)f\|_{L^2}^{\delta(r)}, \quad \forall t \in \mathbb{R},$$
$$\|f\|_{L^r} \leq c_r \left(\frac{\omega}{\sinh(\omega t)}\right)^{\delta(r)}\|f\|_{L^2}^{1-\delta(r)}\|H(t)f\|_{L^2}^{\delta(r)}, \quad \forall t \neq 0.$$

4. They act like derivatives on the nonlinearities $F \in C^1(\mathbb{C}, \mathbb{C})$ satisfying the gauge invariance property $F(z) = G(|z|^2)z$, $\forall z \in \mathbb{C}$, that is,

$$J(t)F(u) = \partial_z F(u)J(t)u - \partial_{\bar{z}}F(u)\overline{J(t)u},$$
$$H(t)F(u) = \partial_z F(u)H(t)v - \partial_{\bar{z}}F(u)\overline{H(t)u}.$$

*Proof.* The first point is easily checked thanks to (1.8). The second assertion is obvious, and imply the last two points. □

*Remark.* One could argue that we consider only isotropic potentials, and not the general form

$$V(x) = \frac{1}{2}\sum_{j=1}^{n}\delta_j \omega_j^2 x_j^2,$$

with $\delta_j \in \{-1, 0, 1\}$, $\omega_j > 0$, not necessarily equal. Strichartz estimates would still be available (locally in time only if some $\delta_j$ is positive), and one could construct operators analogous to $J$ and $H$ that satisfy such properties as those stated in Lemma 2.3. However, the evolution law (2.10) stated below (on which our study highly relies) seems to be bound to isotropic potentials. Finally, the changes of variables we use in Section 5 are also typical of isotropic potentials.

Formally, solutions of (1.2) satisfy the following conservation laws,

(2.8)
$$\text{Mass: } M = \|u(t)\|_{L^2} = cst = \|u_0\|_{L^2},$$
$$\text{Energy: } E = \frac{1}{2}\|\nabla_x u(t)\|_{L^2}^2 - \frac{\omega^2}{2}\|xu(t)\|_{L^2}^2 + \frac{\lambda}{\sigma+1}\|u(t)\|_{L^{2\sigma+2}}^{2\sigma+2} = cst.$$

Notice that even if the nonlinearity is repulsive ($\lambda > 0$), one cannot deduce *a priori* estimates from the conservation of energy. One needs more precise information. Following [Car02b], split the energy into two parts, which are not conserved in general,

$$E_1(t) := \frac{1}{2}\|J(t)u\|_{L^2}^2 + \frac{\lambda}{\sigma+1}\cosh^2(\omega t)\|u(t)\|_{L^{2\sigma+2}}^{2\sigma+2},$$
$$E_2(t) := -\frac{\omega^2}{2}\|H(t)u\|_{L^2}^2 - \frac{\lambda}{\sigma+1}\sinh^2(\omega t)\|u(t)\|_{L^{2\sigma+2}}^{2\sigma+2}.$$

One checks the identity $E_1(t) + E_2(t) \equiv E$.

For $(q,r)$ an admissible pair, define

$$Y_{r,loc}(I) := \{u \in C(I;\Sigma);\ A(t)u \in L^q_{loc}(I;L^r) \cap L^\infty_{loc}(I;L^2),\ \forall A \in \{Id, J, H\}\},$$

$$Y_{loc}(I) := \{u \in C(I;\Sigma);\ A(t)u \in L^q_{loc}(I;L^r),\ \forall A \in \{Id, J, H\}, \forall (q,r) \text{ admissible}\}.$$

When (1.2) is globally well-posed, we will also use

$$Y_r(I) := \{u \in C(I;\Sigma);\ A(t)u \in L^q(I;L^r) \cap L^\infty(I;L^2),\ \forall A \in \{Id, J, H\}\},$$

$$Y(I) := \{u \in C(I;\Sigma);\ A(t)u \in L^q(I;L^r),\ \forall A \in \{Id, J, H\}, \forall (q,r) \text{ admissible}\}.$$

**Proposition 2.4** (Local well-posedness for (1.2)). *Let $\lambda \in \mathbb{R}$, $\sigma, \omega > 0$, with $\sigma < 2/(n-2)$ if $n \geq 3$.*
- *For every $u_0 \in \Sigma$, there exist $t_0 > 0$ independent of $\omega > 0$, and a unique solution $u \in Y_{2\sigma+2,loc}(]-2t_0, 2t_0[)$ to (1.2). Moreover, it belongs to $Y_{loc}(]-2t_0, 2t_0[)$ and there exists $C_0$ depending only on $\lambda$, $n$, $\sigma$ and $\|u_0\|_\Sigma$ such that*

$$(2.9) \qquad \sup_{|t| \leq t_0} \|u(t)\|_{L^2} + \sup_{|t| \leq t_0} \|J(t)u\|_{L^2} + \sup_{|t| \leq t_0} \|H(t)u\|_{L^2} \leq C_0.$$

- *Mass and energy are conserved, that is (2.8) holds. More precisely, $E_1$ and $E_2$ satisfy*

$$(2.10) \qquad \frac{dE_1}{dt} = -\frac{dE_2}{dt} = \frac{\omega\lambda}{2\sigma+2}(2 - n\sigma)\sinh(2\omega t)\|u(t)\|^{2\sigma+2}_{L^{2\sigma+2}}.$$

- *If $u_0^n \to u_0$ in $\Sigma$ and $[t_1, t_2] \subset\,]-2t_0, 2t_0[$, then $u^n \to u$ in $Y([t_1,t_2])$, where $u^n$ solves (1.2) with initial datum $u_0^n$.*

*Proof.* Notice that Duhamel's principle for (1.2) writes

$$(2.11) \qquad u(t) = U_\omega(t)u_0 - i\lambda \int_0^t U_\omega(t-s)\left(|u|^{2\sigma}u\right)(s)ds.$$

The point is to notice that we can reproduce the proof of local existence of solutions to (1.3) in $\Sigma$ (see e.g. [Gin97]). Indeed, Duhamel's principle is similar, we have the same Strichartz inequalities (from Lemma 2.2, with constants independent of $\omega > 0$), and the operators $J$ and $H$ satisfy the same properties as those which are used in the proof of local existence of solutions to (1.3) in $\Sigma$. More precisely, the operators used in the case of (1.3) are $\nabla_x$ and $x + it\nabla_x$; they commute with the linear part of (1.3), act on such nonlinearities as those we consider like derivatives, and $\nabla_x$ provides Gagliardo-Nirenberg inequalities (so does $x + it\nabla_x$, but this point is not used for local existence). Lemma 2.3 shows that $J$ and $H$ satisfy those properties, and the first line of (2.7) provides in particular the same Gagliardo-Nirenberg inequalities as for $\nabla_x$; since $\cosh(x) \geq 1$ for all $x \in \mathbb{R}$,

$$(2.12) \quad \|f\|_{L^r} \leq \frac{c_r}{(\cosh(\omega t))^{\delta(r)}}\|f\|_{L^2}^{1-\delta(r)}\|J(t)f\|_{L^2}^{\delta(r)} \leq c_r\|f\|_{L^2}^{1-\delta(r)}\|J(t)f\|_{L^2}^{\delta(r)}.$$

Moreover, $J(0) = i\nabla_x$ and $H(0) = x$ are independent of $\omega > 0$, so the first point of the proposition follows.

One checks that the identities stated in the second hold for smooth solutions. It follows that they hold for the solutions constructed in the first point, by the same argument as in the case of (1.3) (see e.g. [Caz93], Th. 4.2.8 and Prop. 6.4.2). Similarly, transposing the proof of [Caz93], Th. 4.2.8, yields the last point of the proposition. □

*Remark.* One might think that the above proposition would also hold for (1.6) as well. In that case, we would have some local existence on a time interval independent of $\omega > 0$, which is in contradiction with the result of [Car02b], recalled in the first point of Theorem 1.1. The only difference which makes it impossible to conclude as above, is that in the analog of (2.7), hyperbolic function are replaced by trigonometric functions, that is, the analog of (2.7) is

$$\|f\|_{L^r} \leq \frac{c_r}{(\cos(\omega t))^{\delta(r)}} \|f\|_{L^2}^{1-\delta(r)} \|J(t)f\|_{L^2}^{\delta(r)},$$

$$\|f\|_{L^r} \leq c_r \left(\frac{\omega}{\sin(\omega t)}\right)^{\delta(r)} \|f\|_{L^2}^{1-\delta(r)} \|H(t)f\|_{L^2}^{\delta(r)}.$$

We cannot get rid of the dependence upon $\omega$ as we did in (2.12), which prevents the existence of such a $t_0$ independent of $\omega$.

**Corollary 2.5.** *Let $\lambda \in \mathbb{R}$, $\sigma, \omega > 0$, with $\sigma < 2/(n-2)$ if $n \geq 3$.*
*1. Let $u_0 \in \Sigma$ and $u \in Y_{loc}(I)$ solve (1.2) for some time interval $I$ containing $0$. For any $I \ni t > 0$, the following properties are equivalent:*

- *$\nabla_x u(s)$ is uniformly bounded in $L^2(\mathbb{R}^n)$ for $s \in [0,t]$; $\nabla_x u \in L^\infty([0,t]; L^2)$.*
- *$J(s)u$ or $H(s)u$ is uniformly bounded in $L^2(\mathbb{R}^n)$ for $s \in [0,t]$.*
- *$J(s)u$ and $H(s)u$ are uniformly bounded in $L^2(\mathbb{R}^n)$ for $s \in [0,t]$.*
- *$u(s, \cdot)$ is uniformly bounded in $\Sigma$ for $s \in [0,t]$; $u \in L^\infty([0,t]; \Sigma)$.*

*2. For every $u_0 \in \Sigma$, there exist $T^*(u_0), T_*(u_0) > 0$, and a unique maximal solution $u \in Y_{2\sigma+2,loc}(]-T_*, T^*[)$ to (1.2), which actually belongs to $Y_{loc}(]-T_*, T^*[)$. It is maximal in the sense that if $T^*(u_0) < \infty$, then $\|\nabla_x u(t)\|_{L^2} \to \infty$ as $t \uparrow T^*(u_0)$, and if $T_*(u_0) < \infty$, then $\|\nabla_x u(t)\|_{L^2} \to \infty$ as $t \downarrow -T_*(u_0)$.*

*Proof.* First, notice that the equivalence of the last two properties of the first assertion is a consequence of the formulae (2.3) and (2.4), and of the conservation of mass (2.8).

Assume that $\nabla_x u(s, \cdot)$ is uniformly bounded in $L^2(\mathbb{R}^n)$ for $s \in [0,t]$. Since $u \in Y_{loc}(I)$ solves (1.2), its $L^2$-norm is constant, thus $u(s, \cdot)$ is uniformly bounded in $H^1(\mathbb{R}^n)$ for $s \in [0,t]$. From the Sobolev embedding $H^1(\mathbb{R}^n) \subset L^{2\sigma+2}(\mathbb{R}^n)$, $u(s, \cdot)$ is uniformly bounded in $L^{2\sigma+2}(\mathbb{R}^n)$ for $s \in [0,t]$, and from the conservation of energy (2.8), the first moment of $u$ is uniformly bounded in $L^2(\mathbb{R}^n)$: $u \in L^\infty([0,t]; \Sigma)$.

We now just have to prove that the second and third properties are equivalent, that is, the second implies the third. Assume that $J(s)u$ is uniformly bounded in $L^2(\mathbb{R}^n)$ for $s \in [0,t]$. Then from (2.7), $u(s, \cdot)$ is uniformly bounded in $L^{2\sigma+2}(\mathbb{R}^n)$ for $s \in [0,t]$; $E_1(s)$ is uniformly bounded for $s \in [0,t]$. Since $E_1(s) + E_2(s) \equiv E$, $E_2(s)$ is uniformly bounded for $s \in [0,t]$, which proves that $H(s)u$ is uniformly bounded in $L^2(\mathbb{R}^n)$ for $s \in [0,t]$.

Assume that $H(s)u$ is uniformly bounded in $L^2(\mathbb{R}^n)$ for $s \in [0,t]$. From the first point of Proposition 2.4, $u \in L^\infty([0,t_0]; \Sigma)$ for some $t_0 > 0$. We thus suppose that $t \geq t_0$. From (2.7), $u(s, \cdot)$ is uniformly bounded in $L^{2\sigma+2}(\mathbb{R}^n)$ for $s \in [t_0, t]$, and we can repeat the above argument.

The second assertion follows from the first, and Proposition 2.4. □

We can now state sufficient conditions for the solution of (1.2) to be global. When $\lambda < 0$, let $Q$ denote the unique spherically symmetric solution of (see [Str77], [Kwo89])

(2.13) $$\begin{cases} -\frac{1}{2}\Delta Q + Q = -\lambda |Q|^{4/n} Q, \text{ in } \mathbb{R}^n, \\ Q > 0, \text{ in } \mathbb{R}^n. \end{cases}$$

**Corollary 2.6** (Global existence). *Let $\lambda \in \mathbb{R}$, $\sigma, \omega > 0$, with $\sigma < 2/(n-2)$ if $n \geq 3$, $u_0 \in \Sigma$ and $u \in Y_{loc}(]-T_*, T^*[)$ be the maximal solution given by Corollary 2.5. We have $T_* = T^* = \infty$ in the following cases:*

- *The nonlinearity is repulsive, $\lambda \geq 0$.*

- $\lambda < 0$ and $\sigma < 2/n$.
- $\lambda < 0$, $\sigma = 2/n$ and $\|u_0\|_{L^2} < \|Q\|_{L^2}$.
- $\lambda < 0$, $\sigma > 2/n$ and $\|u_0\|_{H^1}$ is sufficiently small.

In addition, if $\lambda \geq 0$, we have $u \in Y(\mathbb{R})$, and (1.2) is globally well-posed.

*Proof.* We shall prove that under our assumptions, $T^* = \infty$, the proof that $T_* = \infty$ being similar. From Corollary 2.5, it suffices to proves that the $L^2$-norm of $J(t)u$ cannot blow up in finite time.

Assume $\lambda > 0$. If $\sigma \geq 2/n$, then (2.10) implies that for any $t \geq 0$, $E_1(t) \leq E_1(0)$, which yields an *a priori* bound for the $L^2$-norm of $J(t)u$. From Corollary 2.5, this yields $T^* = \infty$. If $\sigma < 2/n$, it follows from (2.10) that for $t \geq 0$,

$$\cosh^2(\omega t)\|u(t)\|_{L^{2\sigma+2}}^{2\sigma+2} \leq \frac{\sigma+1}{\lambda} E_1(0) + \omega\left(1 - \frac{n\sigma}{2}\right)\int_0^t \sinh(2\omega s)\|u(s)\|_{L^{2\sigma+2}}^{2\sigma+2} ds.$$

Gronwall lemma applied to the function defined by the left hand side of the above inequality yields

$$\|u(t)\|_{L^{2\sigma+2}}^{2\sigma+2} \leq \frac{\sigma+1}{\lambda} E_1(0)\big(\cosh(\omega t)\big)^{-n\sigma}.$$

Plugging this estimate in (2.10), we have

(2.14) $$\|J(t)u\|_{L^2}^2 \lesssim e^{(2-n\sigma)\omega t}.$$

This yields $T^* = \infty$.

Now assume $\lambda < 0$. If $\sigma < 2/n$, it follows from (2.10) and (2.7) that for $t \geq 0$,

$$\frac{1}{2}\|J(t)u\|_{L^2}^2 \leq E_1(0) + \frac{|\lambda|}{\sigma+1}\cosh^2(\omega t)\|u(t)\|_{L^{2\sigma+2}}^{2\sigma+2}$$

$$\leq E_1(0) + C\big(\cosh(\omega t)\big)^{2-n\sigma}\|u_0\|_{L^2}^{(2-n)\sigma+2}\|J(t)u\|_{L^2}^{n\sigma}.$$

Since $\sigma < 2/n$, this means that

$$\sup_{t \geq 0}(\cosh(\omega t))^{-2}\|J(t)u\|_{L^2}^2 < \infty,$$

and Corollary 2.5 yields global existence in the future.

If $\sigma = 2/n$, the same argument as above yields

$$\frac{1}{2}\|J(t)u\|_{L^2}^2 \leq E_1(0) + \frac{|\lambda|}{\sigma+1}c_{2+4/n}^{2+4/n}\|u_0\|_{L^2}^{4/n}\|J(t)u\|_{L^2}^2,$$

where $c_{2+4/n}$ is the constant of Gagliardo-Nirenberg inequality mentioned in the third point of Lemma 2.3. Precisely, M. Weinstein ([Wei83]) proved that the best such constant satisfies

$$\frac{|\lambda|}{\sigma+1}c_{2+4/n}^{2+4/n}\|Q\|_{L^2}^{4/n} = \frac{1}{2},$$

where $Q$ is the radial solution of (2.13). Thus if $\|u_0\|_{L^2} < \|Q\|_{L^2}$, we obtain an *a priori* bound for $\|J(t)u\|_{L^2}$, which implies $T^* = \infty$.

Finally, if $\sigma > 2/n$, we have

$$E_1(t) \leq E_1(0) + C\int_0^t \sinh(2\omega s)\|u(s)\|_{L^{2\sigma+2}}^{2\sigma+2} ds$$

$$\leq C\big(\|u_0\|_{H^1}\big) + C\big(\|u_0\|_{L^2}\big)\sup_{0 \leq s \leq t}\|J(s)u\|_{L^2}^{n\sigma}\int_0^t \frac{\sinh(\omega s)}{\cosh^{n\sigma-1}(\omega s)}ds.$$

Therefore,

$$\sup_{0 \leq s \leq t}\|J(s)u\|_{L^2}^2 \leq C\big(\|u_0\|_{H^1}\big) + C\big(\|u_0\|_{L^2}\big)\sup_{0 \leq s \leq t}\|J(s)u\|_{L^2}^{n\sigma},$$

for $C\big(\|u_0\|_{H^1}\big)$ and $C\big(\|u_0\|_{L^2}\big)$ going to zero with their argument. Now we can use the following lemma, whose easy proof is left out.

**Lemma 2.7** (Bootstrap argument). *Let $M = M(t)$ be a nonnegative continuous function on $[0,T]$ such that, for every $t \in [0,T]$,*
$$M(t) \leq a + bM(t)^\theta,$$
*where $a, b > 0$ and $\theta > 1$ are constants such that*
$$a < \left(1 - \frac{1}{\theta}\right) \frac{1}{(\theta b)^{1/(\theta-1)}} \, , \quad M(0) \leq \frac{1}{(\theta b)^{1/(\theta-1)}} \, .$$
*Then, for every $t \in [0,T]$, we have*
$$M(t) \leq \frac{\theta}{\theta - 1} \, a \, \cdot$$

Taking $\|u_0\|_{H^1}$ sufficiently small, we can apply the above lemma and obtain an *a priori* bound for $\|J(t)u\|_{L^2}$.

We now have to prove the last assertion of the corollary, that is, if $\lambda \geq 0$, then
$$(2.15) \qquad A(t)u \in L^q(\mathbb{R}; L^r), \ \forall A \in \{Id, J, H\}, \ \forall (q,r) \text{ admissible}.$$

Let $\lambda \geq 0$ and $A \in \{Id, J, H\}$. In the first part of the proof, we saw that according to the considered case ($\sigma \geq 2/n$, or $\sigma < 2/n$), either $J(t)u \in L^\infty(\mathbb{R}; L^2)$, or it satisfies estimates (2.14). It is easy to check that in either of these two cases, $H(t)u$ satisfies the same estimate as $J(t)u$. Since the second estimate is weaker than the first one, it suffices to prove that it yields (2.15). We will use the following algebraic lemma,

**Lemma 2.8.** *Let $r = s = 2\sigma + 2$, and $q$ be such that the pair $(q,r)$ is admissible. Define $k$ by*
$$k = \frac{2\sigma(2\sigma + 2)}{2 - (n-2)\sigma}.$$
*Then $k$ is finite, and the following algebraic identities hold,*
$$\begin{cases} \dfrac{1}{r'} = \dfrac{1}{r} + \dfrac{2\sigma}{s} \, , \\ \dfrac{1}{q'} = \dfrac{1}{q} + \dfrac{2\sigma}{k} \, . \end{cases}$$

Let $q, r, k$ and $s$ be as in the above lemma. From (2.14) and the conservation of mass, (2.7) yields
$$(2.16) \qquad \|u\|_{L^k([T,\infty[;L^s)} \leq C \left\|e^{-n\sigma\delta(s)t/2}\right\|_{L^k([T,\infty[)} \leq Ce^{-n\sigma\delta(s)T/2} \, .$$

To prove that $A(t)u \in L^q([0,\infty[; L^r)$, write Duhamel's principle with time origin equal to $T$, to be fixed later,
$$u(t) = U_\omega(t-T)u(T) - i\lambda \int_T^t U_\omega(t-s)\left(|u|^{2\sigma}u\right)(s)ds.$$

Applying the operator $A$ yields,
$$A(t)u = U_\omega(t-T)A(T)u - i\lambda \int_T^t U_\omega(t-s)A(s)\left(|u|^{2\sigma}u\right)ds,$$

and from Strichartz inequalities and Lemma 2.8, for any $S > T$,
$$\|A(t)u\|_{L^q([T,S];L^r)} \leq C_r \|A(T)u\|_{L^2} + C_{r,r} \left\|A(t)\left(|u|^{2\sigma}u\right)\right\|_{L^{q'}([T,S];L^{r'})}$$
$$\leq C_r \|A(T)u\|_{L^2} + \underline{C}\|u\|_{L^k([T,\infty[;L^s)}^{2\sigma} \|A(t)\|_{L^q([T,S];L^r)},$$

where $\underline{C}$ does not depend on $T, S$. From (2.16), choosing $T$ sufficiently large, the second term of the right hand side can be absorbed by the left hand side, and
$$\|A(t)u\|_{L^q([T,S];L^r)} \leq 2C_r \|A(T)u\|_{L^2}.$$

Since $S > T$ is arbitrary, this implies $A(t)u \in L^q(\mathbb{R}_+; L^r)$. Similarly, $A(t)u \in L^q(\mathbb{R}; L^r)$; (2.15) is proven for the admissible pair $(q, r)$ such that $r = 2\sigma + 2$. Strichartz inequality (2.2) then yields (2.15) for any admissible pair. Indeed, if $(q_1, r_1)$ is admissible,

$$\|A(t)u\|_{L^{q_1}([0,S];L^{r_1})} \leq C_{r_1}\|u_0\|_\Sigma + |\lambda|C_{r,1,r}\left\|A(t)\left(|u|^{2\sigma}u\right)\right\|_{L^{q'}([0,S];L^{r'})}$$
$$\leq C + C\|u\|^{2\sigma}_{L^k(\mathbb{R};L^s)}\|A(t)\|_{L^q(\mathbb{R};L^r)}.$$

This completes the proof of Corollary 2.6. □

## 3. Scattering theory

In this section, we prove that the influence of the nonlinear term in (1.2) is negligible as time becomes large (at least if $\lambda > 0$), without the usual restriction on the power of the nonlinearity encountered for scattering theory associated to (1.3). We first prove the existence of wave operators, then their asymptotic completeness.

**Proposition 3.1** (Existence of wave operators). *In either of the cases considered in Corollary 2.6, the following holds.*

- *For every $u_- \in \Sigma$, there exists a unique $u_0 \in \Sigma$ such that the maximal solution $u \in C(\mathbb{R}; \Sigma)$ of (1.2) satisfies*

$$\left\|U_\omega(-t)u(t) - u_-\right\|_\Sigma \xrightarrow[t \to -\infty]{} 0.$$

- *For every $u_+ \in \Sigma$, there exists a unique $u_0 \in \Sigma$ such that the maximal solution $u \in C(\mathbb{R}; \Sigma)$ of (1.2) satisfies*

$$\left\|U_\omega(-t)u(t) - u_+\right\|_\Sigma \xrightarrow[t \to +\infty]{} 0.$$

*Proof.* We prove the first point, the proof of the second being similar. We solve the following equation by a fixed point argument,

$$(3.1) \qquad u(t) = U_\omega(t)u_- - i\lambda \int_{-\infty}^t U_\omega(t-s)\bigl(|u|^{2\sigma}u\bigr)(s)ds.$$

Define $F(u)(t)$ as the right hand side of (3.1), and let $R := \|u_-\|_\Sigma$. Recall that as stated in Lemma 2.8, $(q, r)$ is the admissible pair such that $r = 2\sigma + 2$.

We first prove that there exists $T > 0$ such that the set

$$X_T := \{u \in Y_{2\sigma+2}(]-\infty, -T]);\ \|A(t)u\|_{L^2} \leq 2R\ ,\ \forall t \leq -T,\ A \in \{Id, J, H\},$$
$$\|A(t)u\|_{L^q(]-\infty,-T];L^r)} \leq 2C_{2\sigma}R,\ \forall A \in \{Id, J, H\}\}$$

is stable under the map $F$, where $C_{2\sigma}$ is the constant in Strichartz inequality (2.1). We then prove that up to choosing $T$ even larger, $F$ is a contraction on $X_T$, which is equipped with the norm

$$\|u\|_{X_T} := \sum_{A \in \{Id, J, H\}} \bigl(\|A(t)u\|_{L^\infty(]-\infty,-T];L^2)} + \|A(t)u\|_{L^q(]-\infty,-T];L^r)}\bigr).$$

For any pair $(a, b)$, we use the notation

$$\|f\|_{L^a_T(L^b)} = \|f\|_{L^a(]-\infty,-T];L^b)}.$$

Let $u \in X_T$, and $A \in \{Id, J, H\}$. From Lemma 2.2, Lemma 2.3 and Lemma 2.8,

$$\|A(t)F(u)\|_{L^\infty_T(L^2)} \leq \|u_-\|_\Sigma + C_{2,2\sigma+2}|\lambda|\left\|A(t)\left(|u|^{2\sigma}u\right)\right\|_{L^{q'}_T(L^{r'})}$$
$$\leq R + C\left\||u|^{2\sigma}A(t)u\right\|_{L^{q'}_T(L^{r'})}$$
$$\leq R + C\|u\|^{2\sigma}_{L^k_T(L^s)}\|A(t)u\|_{L^q_T(L^r)}.$$

From (2.7) and Lemma 2.8,

$$\|u\|_{L_T^k(L^s)} \leq C_k R \left\|\left(\frac{1}{\cosh(\omega t)}\right)^{\delta(s)}\right\|_{L^k(]-\infty,-T])} \leq C(\omega,\sigma) R e^{-\omega\delta(s)T}.$$

It follows,

(3.2) $$\|A(t)F(u)\|_{L_T^\infty(L^2)} \leq R + CR^{2\sigma+1}e^{-2\sigma\omega\delta(s)T}.$$

Use Lemma 2.2, Lemma 2.3 and Lemma 2.8 again to obtain

$$\|A(t)F(u)\|_{L_T^q(L^r)} \leq C_{2\sigma}R + C \|u\|_{L_T^k(L^s)}^{2\sigma} \|A(t)u\|_{L_T^q(L^r)}$$
$$\leq C_{2\sigma}R + CR^{2\sigma+1}e^{-2\sigma\omega\delta(s)T}.$$

It is now clear that if $T$ is sufficiently large, then $X_T$ is stable under $F$.

To complete the proof of the proposition, following the argument used in [Kat87], it is enough to prove contraction for large $T$ in the weaker metric $L^q(]-\infty,-T];L^r)$. From Lemma 2.2, Lemma 2.3 and Lemma 2.8, we have

(3.3) $$\|F(u_2) - F(u_1)\|_{L_T^q(L^r)} \leq C \left\|\left(|u_2|^{2\sigma}u_2 - |u_1|^{2\sigma}u_1\right)\right\|_{L_T^{q'}(L^{r'})}$$
$$\leq C \left(\|u_1\|_{L_T^k(L^s)}^{2\sigma} + \|u_2\|_{L_T^k(L^s)}^{2\sigma}\right) \|u_2 - u_1\|_{L_T^q(L^r)}.$$

As above, we have the estimate

$$\|u_1\|_{L_T^k(L^s)}^{2\sigma} + \|u_2\|_{L_T^k(L^s)}^{2\sigma} \leq CR^{2\sigma}e^{-2\sigma\omega\delta(s)T}.$$

Therefore, contraction follows for $T$ sufficiently large.

From Corollary 2.6, the solution $u$ we obtain by this fixed point argument is defined not only on $]-\infty,-T]$ for $T$ large, but globally. Proposition 3.1 then follows from Corollaries 2.5 and 2.6. □

*Remark.* The fact that we limit ourselves to the cases considered in Corollary 2.6 in the above proposition is needed only to ensure the solution $u$ we construct is defined up to time $t=0$. To solve (3.1) in the neighborhood of $-\infty$, we used only the assumptions of Proposition 2.4.

**Proposition 3.2** (Asymptotic completeness). *Let $\lambda \geq 0$, $\sigma > 0$, with $\sigma < 2/(n-2)$ if $n \geq 3$.*

- *For every $u_0 \in \Sigma$, there exists a unique $u_- \in \Sigma$ such that the maximal solution $u \in C(\mathbb{R};\Sigma)$ of (1.2) satisfies*

$$\|U_\omega(-t)u(t) - u_-\|_\Sigma \xrightarrow[t\to-\infty]{} 0.$$

- *For every $u_0 \in \Sigma$, there exists a unique $u_+ \in \Sigma$ such that the maximal solution $u \in C(\mathbb{R};\Sigma)$ of (1.2) satisfies*

$$\|U_\omega(-t)u(t) - u_+\|_\Sigma \xrightarrow[t\to+\infty]{} 0.$$

*Proof.* We prove the second point, the proof of the first being similar. Let $u_0 \in \Sigma$.

Since $\Sigma$ is a Hilbert space, it is enough to prove that the family $\left(U_\omega(-t)u(t)\right)_{t\geq 0}$ is a Cauchy sequence as $t$ goes to $+\infty$. From Duhamel's principle (2.11), we have

$$U_\omega(-t)u(t) = u_0 - i\lambda \int_0^t U_\omega(-s) \left(|u|^{2\sigma}u\right)(s)ds.$$

Let $B \in \{Id, \nabla_x, x\}$, and $A \in \{Id, J, H\}$ be its counterpart given by the commutation property (2.5). We have

$$B\left(U_\omega(-t)u(t)\right) = Bu_0 - i\lambda \int_0^t U_\omega(-s)A(s)\left(|u|^{2\sigma}u\right)ds.$$

Let $t_2 \geq t_1 > 0$. From Strichartz inequality (2.2),

$$\left\|B\big(U_\omega(-t_2)u(t_2) - U_\omega(-t_1)u(t_1)\big)\right\|_{L^2} \leq$$

$$\leq C \left\| \int_{t_1}^t U_\omega(-s) A(s)\left(|u|^{2\sigma} u\right) ds \right\|_{L^\infty([t_1,t_2];L^2)}$$

$$\leq C \left\| A(t)\left(|u|^{2\sigma} u\right) \right\|_{L^{q'}([t_1,t_2];L^{r'})}.$$

From Lemma 2.3 and Lemma 2.8, this yields

$$\left\|B\big(U_\omega(-t_2)u(t_2) - U_\omega(-t_1)u(t_1)\big)\right\|_{L^2} \leq C \|u\|_{L^k([t_1,t_2];L^s)}^{2\sigma} \|A(t)u\|_{L^q([t_1,t_2];L^r)}$$

$$\leq C \|u\|_{L^k([t_1,t_2];L^s)}^{2\sigma} \|A(t)u\|_{L^q([t_1,t_2];L^r)}.$$

We saw in the proof of Corollary 2.6 that $u \in L^k(\mathbb{R}; L^s)$ and $A(t)u \in L^q(\mathbb{R}; L^r)$, which implies that $\big(B(U_\omega(-t)u(t))\big)_{t>0}$ is a Cauchy sequence in $L^2$, and completes the proof of the proposition. $\square$

Propositions 3.1 and 3.2 imply Theorem 1.3, and even more, since we do not necessarily assume that the nonlinearity is defocusing.

## 4. Blow-up in finite time

In Corollary 2.6, we proved that if $\lambda < 0$ and $\sigma \geq 2/n$, the solution of (1.2) is global provided that the initial datum $u_0$ is small. When $u_0$ is not small, we show that finite time blow up may occur, as in the case of (1.3). However, the sufficient condition we state below is stronger than its counterpart (1.4) for (1.3); in some sense, blow-up in finite time is less likely to occur for (1.2) than for (1.3).

**Proposition 4.1.** *Let $u_0 \in \Sigma$, $\lambda < 0$, $\sigma \geq 2/n$, with $\sigma < 2/(n-2)$ if $n \geq 3$. If $u_0$ satisfies*

$$\frac{1}{2}\|\nabla u_0\|_{L^2}^2 + \frac{\lambda}{\sigma+1}\|u_0\|_{L^{2\sigma+2}}^{2\sigma+2} < -\frac{\omega^2}{2}\|xu_0\|_{L^2}^2,$$

*then the solution $u$ to (1.2) blows up in finite time, in the future or in the past. More precisely,*

- *If $\operatorname{Im} \int \overline{u_0} x \cdot \nabla u_0 \leq 0$, then $T^* < \infty$, that is, $u$ blows up in the future.*
- *If $\operatorname{Im} \int \overline{u_0} x \cdot \nabla u_0 \geq 0$, then $T_* < \infty$, that is, $u$ blows up in the past.*

*If moreover*

$$\frac{1}{2}\|\nabla u_0\|_{L^2}^2 + \frac{\lambda}{\sigma+1}\|u_0\|_{L^{2\sigma+2}}^{2\sigma+2} < -\frac{\omega^2}{2}\|xu_0\|_{L^2}^2 - \omega \left|\operatorname{Im}\int \overline{u_0} x \cdot \nabla_x u_0 \right|,$$

*then $u$ blows up in the past and in the future.*

*Proof.* We follow the method of Zakharov-Glassey. Denote $y(t) := \|xu(t)\|_{L^2}^2$. We show that $y(t)$ satisfies a second order ordinary differential equation, from which the proposition follows.

**Step 1. Formal computations.** Differentiating $y(t)$ and using Eq. (1.2) yields

$$\dot{y}(t) = 2 \operatorname{Im} \int \overline{u}(t,x) x \cdot \nabla_x u(t,x) dx.$$

Expanding $\|J(t)u\|_{L^2}^2$, we have

$$\|J(t)u\|_{L^2}^2 = \omega^2 \sinh^2(\omega t) y(t) + \cosh^2(\omega t)\|\nabla_x u(t)\|_{L^2}^2$$

$$- \omega \sinh(2\omega t) \operatorname{Im} \int \overline{u}(t,x) x \cdot \nabla_x u(t,x) dx,$$

and from the conservation of energy (2.8),
$$E_1(t) = \frac{\omega^2}{2}\sinh^2(\omega t) y(t) + \cosh^2(\omega t)\left(E + \frac{\omega^2}{2}y(t)\right)$$
$$- \frac{\omega}{2}\sinh(2\omega t)\,\mathrm{Im}\int \overline{u}(t,x) x \cdot \nabla_x u(t,x)dx.$$

Using the evolution law (2.10) and the above formula for $\dot{y}(t)$ yields,
$$\frac{d}{dt}\mathrm{Im}\int \overline{u}(t,x)x\cdot\nabla_x u(t,x)dx = 2\omega^2 y(t) + 2E - \frac{\lambda}{\sigma+1}(2-n\sigma)\|u(t)\|_{L^{2\sigma+2}}^{2\sigma+2}.$$

It follows,

(4.1) $$\ddot{y}(t) = 4\omega^2 y(t) + 4E - \frac{2\lambda}{\sigma+1}(2-n\sigma)\|u(t)\|_{L^{2\sigma+2}}^{2\sigma+2}.$$

**Step 2. Justification.** One has to know that $y \in C^1(]-T_*, T^*[)$, the rest of the computations follow from (2.10). The argument is classical (it consists in considering $y^\varepsilon(t) := \|e^{-\varepsilon|x|^2} xu(t)\|_{L^2}^2$ and eventually letting $\varepsilon$ go to zero), and we refer to [Caz93], Section 6.4, for more details, as we did for the proof of (2.10).

**Step 3. Conclusion.** From classical ordinary differential equations' methods, the solution of (4.1) is given by the formula
$$y(t) = y(0)\cosh(2\omega t) + \dot{y}(0)\frac{\sinh(2\omega t)}{2\omega} + \int_0^t \frac{\sinh(2\omega(t-s))}{2\omega}f(s)ds,$$
where $f(t) = 4E - \frac{2\lambda}{\sigma+1}(2-n\sigma)\|u(t)\|_{L^{2\sigma+2}}^{2\sigma+2}$. Since $\lambda < 0$ and $\sigma \geq 2/n$,
$$y(t) \leq y(0)\cosh(2\omega t) + \dot{y}(0)\frac{\sinh(2\omega t)}{2\omega} + \int_0^t \frac{\sinh(2\omega(t-s))}{2\omega}4E\,ds$$
$$\leq y(0)\cosh^2(\omega t) + \frac{2\sinh^2(\omega t)}{\omega^2}\left(\frac{1}{2}\|\nabla u_0\|_{L^2}^2 + \frac{\lambda}{\sigma+1}\|u_0\|_{L^{2\sigma+2}}^{2\sigma+2}\right)$$
$$+ \dot{y}(0)\frac{\sinh(2\omega t)}{2\omega}.$$

Assume $\dot{y}(0) \leq 0$. Then for positive times, the above estimate implies
$$y(t) \leq y(0)\cosh^2(\omega t) + \frac{2\sinh^2(\omega t)}{\omega^2}\left(\frac{1}{2}\|\nabla u_0\|_{L^2}^2 + \frac{\lambda}{\sigma+1}\|u_0\|_{L^{2\sigma+2}}^{2\sigma+2}\right)$$
$$\leq \cosh^2(\omega t)\left(y(0) + \frac{2\tanh^2(\omega t)}{\omega^2}\left(\frac{1}{2}\|\nabla u_0\|_{L^2}^2 + \frac{\lambda}{\sigma+1}\|u_0\|_{L^{2\sigma+2}}^{2\sigma+2}\right)\right).$$

Since $\tanh(\mathbb{R}_+) = [0,1[$, it follows that if
$$\frac{1}{2}\|\nabla u_0\|_{L^2}^2 + \frac{\lambda}{\sigma+1}\|u_0\|_{L^{2\sigma+2}}^{2\sigma+2} < -\frac{\omega^2}{2}\|xu_0\|_{L^2}^2,$$
then if we suppose $T^* = \infty$, $y(t)$ becomes negative for possibly large $t$. This is absurd, therefore $T^*$ is finite. Similarly, if $\dot{y}(0) \geq 0$, then $T_*$ is finite.

Finally, if
$$\frac{1}{2}\|\nabla u_0\|_{L^2}^2 + \frac{\lambda}{\sigma+1}\|u_0\|_{L^{2\sigma+2}}^{2\sigma+2} < -\frac{\omega^2}{2}\|xu_0\|_{L^2}^2 - \omega\left|\mathrm{Im}\int \overline{u_0}x\cdot\nabla_x u_0\right|,$$
then the same argument shows that $T_*$ and $T^*$ are finite. $\square$

We now prove that indeed, blow-up in finite time is less likely to occur for (1.2) than for (1.3). For a fixed $u_0 \in \Sigma$, the blow-up sufficient conditions stated in Proposition 4.1 become empty when $\omega$ is large. We prove that for a fixed initial datum $u_0$, taking $\omega$ sufficiently large ensures the global existence of $u$.

**Proposition 4.2.** *Let $u_0 \in \Sigma$, $\lambda < 0$, $\sigma > 2/n$, with $\sigma < 2/(n-2)$ if $n \geq 3$. There exists $\omega_1 > 0$ such that for any $\omega \geq \omega_1$, the solution $u$ to (1.2) is global, and $u \in Y(\mathbb{R})$.*

*Proof.* From Proposition 2.4, there exist $t_0 > 0$ and $C_0$ independent of $\omega > 0$ such that (1.2) has a unique solution $u \in Y(]-2t_0, 2t_0[)$, which satisfies in addition the estimate (2.9), that is

$$\sup_{|t| \leq t_0} \|u(t)\|_{L^2} + \sup_{|t| \leq t_0} \|J(t)u\|_{L^2} + \sup_{|t| \leq t_0} \|H(t)u\|_{L^2} \leq C_0 .$$

The idea is to mimic the proof of the fourth case in Corollary 2.6, by replacing the smallness assumption by the property $\omega \gg 1$.

Integrate the evolution law (2.10) between time $t_0$ and time $t > t_0$:

$$E_1(t) - E_1(t_0) = \frac{\omega \lambda}{2\sigma + 2}(2 - n\sigma) \int_{t_0}^{t} \sinh(2\omega s) \|u(s)\|_{L^{2\sigma+2}}^{2\sigma+2} ds.$$

From Proposition 2.4 and the fact that the nonlinearity we consider is focusing, we have

$$E_1(t_0) \leq \frac{1}{2}\|J(t_0)u\|_{L^2}^2 \leq \frac{1}{2}C_0^2,$$

where $C_0$ does not depend on $\omega$.

Using modified Gagliardo-Nirenberg inequalities (2.7), we have

$$E_1(t) \leq \frac{1}{2}C_0^2 + C(\lambda, \sigma)\omega \int_{t_0}^{t} \frac{\sinh(2\omega s)}{(\cosh(\omega s))^{n\sigma}} \|u(s)\|_{L^2}^{2+(2-n)\sigma} \|J(s)u\|_{L^2}^{n\sigma} ds$$

$$\leq \frac{1}{2}C_0^2 + C'(\lambda, \sigma)\|u_0\|_{L^2}^{2+(2-n)\sigma} \sup_{t_0 \leq s \leq t} \|J(s)u\|_{L^2}^{n\sigma} (\cosh(\omega t_0))^{2-n\sigma}.$$

From the definition of $E_1$, this yields,

$$\frac{1}{2}\|J(t)u\|_{L^2}^2 \leq \frac{1}{2}C_0^2 + C(\cosh(\omega t))^2 \|u(t)\|_{L^{2\sigma+2}}^{2\sigma+2}$$
$$+ C \sup_{t_0 \leq s \leq t} \|J(s)u\|_{L^2}^{n\sigma} (\cosh(\omega t_0))^{2-n\sigma}$$
$$\leq \frac{1}{2}C_0^2 + C(\cosh(\omega t))^{2-n\sigma} \|J(t)u\|_{L^2}^{n\sigma}$$
$$+ C \sup_{t_0 \leq s \leq t} \|J(s)u\|_{L^2}^{n\sigma} (\cosh(\omega t_0))^{2-n\sigma}$$
$$\leq \frac{1}{2}C_0^2 + C(\cosh(\omega t_0))^{2-n\sigma} \|J(t)u\|_{L^2}^{n\sigma}$$
$$+ C \sup_{t_0 \leq s \leq t} \|J(s)u\|_{L^2}^{n\sigma} (\cosh(\omega t_0))^{2-n\sigma} ,$$

where the above constants do not depend on $\omega$. We finally obtain,

$$\sup_{t_0 \leq s \leq t} \|J(s)u\|_{L^2}^2 \leq C_0^2 + C(\cosh(\omega t_0))^{2-n\sigma} \sup_{t_0 \leq s \leq t} \|J(s)u\|_{L^2}^{n\sigma} ,$$

which can also be written

$$\sup_{t_0 \leq s \leq t} \|J(s)u\|_{L^2}^2 \leq C_0^2 + f(\omega) \left( \sup_{t_0 \leq s \leq t} \|J(s)u\|_{L^2}^2 \right)^{n\sigma/2} ,$$

where $C_0$ does not depend on $\omega$, and $f(\omega) \to 0$ as $\omega \to +\infty$, because $n\sigma > 2$. Lemma 2.7 shows that for $\omega$ sufficiently large, $J(t)u$ is uniformly bounded in $L^2$ for $t \geq 0$. Corollary 2.5 then implies $u \in Y_{loc}(\mathbb{R}_+)$, and one can repeat the end of the proof of Corollary 2.6 to deduce that $u \in Y(\mathbb{R}_+)$.

Proving $u \in Y(\mathbb{R}_-)$ is similar, and we leave out this part. $\square$

The proof of Theorem 1.1 is not complete yet, for in the above proposition, we assumed only $\sigma > 2/n$, while in Theorem 1.1, we assumed $\sigma \geq 2/n$. The remaining case $\sigma = 2/n$ is treated in the next section.

## 5. The particular case $\sigma = 2/n$

Let $\lambda \in \mathbb{R}$, $u_0 \in \Sigma$. Let $v$ solve (1.3) with a critical power, that is,

$$
(5.1) \quad \begin{cases} i\partial_t v + \dfrac{1}{2}\Delta v = \lambda |v|^{4/n} v, & (t,x) \in \mathbb{R} \times \mathbb{R}^n, \\ v_{|t=0} = u_0. \end{cases}
$$

Let $\omega > 0$. In [Car02a], we noticed that if $u^+$ is defined by

$$
(5.2) \quad u^+(t,x) = \frac{1}{(\cos(\omega t))^{n/2}} e^{-i\frac{\omega}{2}|x|^2 \tan(\omega t)} v\left(\frac{\tan(\omega t)}{\omega}, \frac{x}{\cos(\omega t)}\right),
$$

then $u^+$ solves (1.6) with $\sigma = 2/n$. We also proved that $v$ blows up at time $T > 0$ if and only if $u^+$ blows up at time $\arctan(\omega T)/\omega$. The first point of Theorem 1.2 is therefore a reminder of a result stated in [Car02a].

As noticed in the introduction, replacing $\omega$ by $\pm i\omega$ formally turns (1.6) into (1.2). Following this idea again, define

$$
(5.3) \quad u^-(t,x) = \frac{1}{(\cosh(\omega t))^{n/2}} e^{i\frac{\omega}{2}|x|^2 \tanh(\omega t)} v\left(\frac{\tanh(\omega t)}{\omega}, \frac{x}{\cosh(\omega t)}\right).
$$

Then from Proposition 2.4, $u^-$ is the solution of

$$
(5.4) \quad \begin{cases} i\partial_t u^- + \dfrac{1}{2}\Delta u^- = -\omega^2 \dfrac{|x|^2}{2} u^- + \lambda |u^-|^{4/n} u^-, & (t,x) \in \mathbb{R} \times \mathbb{R}^n, \\ u^-_{|t=0} = u_0. \end{cases}
$$

Now assume that $\lambda < 0$, and that $v$ blows up at some finite time $T > 0$.

From the factorization (2.6), it is easy to see that

$$
(5.5) \quad \|J(t)u^-\|_{L^2} = \left\|\nabla_x v\left(\frac{\tanh(\omega t)}{\omega}\right)\right\|_{L^2}.
$$

Since $\tanh(\mathbb{R}_+) = [0,1[$, if $\omega \geq 1/T$, then the function of the right hand side of (5.3) does not "see" the time $T$, and from Corollary 2.5, $u^-$ does not blow up in finite time.

If $\omega < 1/T$, then (5.5) and Corollary 2.5 show that $u^-$ blows up at time

$$
T_\omega = \frac{\arg\tanh(\omega T)}{\omega},
$$

which completes the proof of Theorem 1.2.

We can go further in the analysis of the influence of the parameter $\omega$.

**Proposition 5.1.** *Let $u_0 \in \Sigma$, $\lambda < 0$. For $\omega \geq 0$, denote $u^\omega$ the solution of the initial value problem*

$$
(5.6) \quad \begin{cases} i\partial_t u^\omega + \dfrac{1}{2}\Delta u^\omega = -\omega^2 \dfrac{|x|^2}{2} u^\omega + \lambda |u^\omega|^{4/n} u^\omega, & (t,x) \in \mathbb{R} \times \mathbb{R}^n, \\ u^\omega_{|t=0} = u_0. \end{cases}
$$

*Let $\omega_* \geq 0$.*
• *If $u^{\omega_*}$ is defined globally, $u^{\omega_*} \in Y_{loc}(\mathbb{R})$, then for every $\omega \geq \omega_*$, $u^\omega$ is also defined globally, $u^\omega \in Y_{loc}(\mathbb{R})$.*
• *Suppose that there exists $T_* > 0$ such that $u^{\omega_*}$ blows up at time $T_*$. Then for every $0 \leq \omega < \omega_*/\tanh(\omega_* T_*)$, $u^\omega$ blows up in finite time, and for every $\omega \geq \omega_*/\tanh(\omega_* T_*)$, $u^\omega$ is defined globally, $u^\omega \in Y_{loc}(\mathbb{R})$.*

*Proof.* Let $v$ be the solution of (5.1). If $v$ is defined globally in $\Sigma$, then so is $u^\omega$ for any $\omega \geq 0$. If $v$ blows up in finite time $T_0 > 0$ while $u^{\omega_*}$ is defined globally, then Theorem 1.2 implies $\omega_* \geq 1/T_0$. Using Theorem 1.2 again, $u^\omega$ is defined globally for any $\omega \geq \omega_* \geq 1/T_0$.

Now assume $u^{\omega_*}$ blows up at time $T_* > 0$. From Theorem 1.2, $v$ blows up at time $T_0$, with

$$T_0 = \frac{\tanh(\omega_* T_*)}{\omega_*} \ .$$

The last point of Theorem 1.2 implies that if $\omega \geq 1/T_0$, then $u^\omega$ is defined globally. Similarly, if $0 \leq \omega < \omega_*/\tanh(\omega_* T_*)$, then $u^\omega$ blows up at time

$$T_\omega = \frac{1}{\omega} \arg\tanh\left(\frac{\omega \tanh(\omega_* T_*)}{\omega_*}\right),$$

which completes the proof of the proposition. $\square$

Assume $\lambda < 0$, and that $v$ blows up in finite time. Theorem 1.2 and Proposition 5.1 show that there is a critical value for the parameter $\omega$, which is the inverse of the blow-up time for $v$. What happens for that critical value? We can answer this question in the case where the mass of the initial datum is critical. We saw in Corollary 2.6 that if $\|u_0\|_{L^2} < \|Q\|_{L^2}$, where $Q$ is the spherically symmetric solution of (2.13), then the solution to (5.6) is global for any $\omega \geq 0$. If $\|u_0\|_{L^2} = \|Q\|_{L^2}$, then blow-up in finite time may occur. This phenomenon was studied very precisely by Merle in the case of (5.1).

**Theorem 5.2** ([Mer93], Th. 1). *Let $\lambda < 0$, $u_0 \in H^1(\mathbb{R}^n)$, and assume that the solution $v$ of (5.1) blows up in finite time $T > 0$. Moreover, assume that $\|u_0\|_{L^2} = \|Q\|_{L^2}$, where $Q$ is defined by (2.13). Then there exist $\theta \in \mathbb{R}$, $\delta > 0$, $x_0, x_1 \in \mathbb{R}^n$ such that*

$$u_0(x) = \left(\frac{\delta}{T}\right)^{n/2} e^{i\theta - i|x-x_1|^2/2T + i\delta^2/T} Q\left(\delta\left(\frac{x-x_1}{T} - x_0\right)\right),$$

*and for $t < T$,*

$$v(t,x) = \left(\frac{\delta}{T-t}\right)^{n/2} e^{i\theta - i|x-x_1|^2/2(T-t) + i\delta^2/(T-t)} Q\left(\delta\left(\frac{x-x_1}{T-t} - x_0\right)\right).$$

We use this result only to understand the role of $\omega$ to prevent blow-up when the mass is critical, but other applications are possible (see [Car02a] for the case of a confining harmonic potential). With the above theorem and the change of variable (5.3), the following result is straightforward.

**Corollary 5.3.** *Let $\lambda < 0$ and $T > 0$. Assume that $u_0$ is given by*

$$u_0(x) = \frac{1}{T^{n/2}} e^{-i|x|^2/2T + i/T} Q\left(\frac{x}{T}\right).$$

*For $\omega \geq 0$, denote $u^\omega$ the solution of (5.6).*
- *If $0 \leq \omega < 1/T$, then $u^\omega$ blows up at time $\arg\tanh(\omega T)/\omega$, with the profile $Q$.*
- *If $\omega > 1/T$, then $u^\omega$ is defined globally, with exponential decay, $u^\omega \in Y(\mathbb{R})$.*
- *If $\omega = 1/T$, then $u^\omega$ is defined globally, with only $u^\omega \in Y_{loc}(\mathbb{R})$. More precisely,*

$$u^{1/T}(t,x) = \left(\omega e^{\omega t}\right)^{n/2} Q\left(\omega x e^{\omega t}\right) e^{-i\omega|x|^2/2 + i\omega(e^{2\omega t}+1)/2}$$

$$= \left(\frac{e^{t/T}}{T}\right)^{n/2} Q\left(\frac{x e^{t/T}}{T}\right) e^{-i|x|^2/2T + i(e^{2t/T}+1)/2T} \ .$$

The critical value $\omega = 1/T$ thus leads to a global solution (we already knew that, from Theorem 1.2), which may have exponential growth (and does, in the particular case $\|u_0\|_{L^2} = \|Q\|_{L^2}$).

# References

[Bar84] J. E. Barab, *Nonexistence of asymptotically free solutions for nonlinear Schrödinger equation*, J. Math. Phys. **25** (1984), 3270–3273.

[Car01] R. Carles, *Geometric optics and long range scattering for one-dimensional nonlinear Schrödinger equations*, Comm. Math. Phys. **220** (2001), no. 1, 41–67.

[Car02a] ———, *Critical nonlinear Schrödinger equations with and without harmonic potential*, Math. Mod. Meth. Appl. Sci. **12** (2002), no. 10, 1513–1523.

[Car02b] ———, *Remarks on nonlinear Schrödinger equations with harmonic potential*, Ann. Henri Poincaré **3** (2002), no. 4, 757–772.

[Caz93] T. Cazenave, *An introduction to nonlinear Schrödinger equations*, Text. Met. Mat., vol. 26, Univ. Fed. Rio de Jan., 1993.

[CFG] R. Carles, C. Fermanian, and I. Gallagher, *On the role of quadratic oscillations in nonlinear Schrödinger equations*, `arXiv:math.AP/0212171`, to appear in J. Funct. Anal.

[CW92] T. Cazenave and F. Weissler, *Rapidly decaying solutions of the nonlinear Schrödinger equation*, Comm. Math. Phys. **147** (1992), 75–100.

[DDM02] Arnaud Debussche and Laurent Di Menza, *Numerical simulation of focusing stochastic nonlinear Schrödinger equations*, Phys. D **162** (2002), no. 3-4, 131–154.

[DS63] Nelson Dunford and Jacob T. Schwartz, *Linear operators. Part II: Spectral theory. Self adjoint operators in Hilbert space*, With the assistance of William G. Bade and Robert G. Bartle, Interscience Publishers John Wiley & Sons New York-London, 1963.

[FH65] R. P. Feynman and A.R. Hibbs, *Quantum mechanics and path integrals (International Series in Pure and Applied Physics)*, Maidenhead, Berksh.: McGraw-Hill Publishing Company, Ltd., 365 p., 1965.

[Gin97] J. Ginibre, *An introduction to nonlinear Schrödinger equations*, Nonlinear waves (Sapporo, 1995) (R. Agemi, Y. Giga, and T. Ozawa, eds.), GAKUTO International Series, Math. Sciences and Appl., Gakkōtosho, Tokyo, 1997, pp. 85–133.

[Gla77] R. T. Glassey, *On the blowing up of solutions to the Cauchy problem for nonlinear Schrödinger equations*, J. Math. Phys. **18** (1977), 1794–1797.

[HN98] N. Hayashi and P. Naumkin, *Asymptotics for large time of solutions to the nonlinear Schrödinger and Hartree equations*, Amer. J. Math. **120** (1998), 369–389.

[Kat87] T. Kato, *Nonlinear Schrödinger equations*, Ann. IHP (Phys. Théor.) **46** (1987), 113–129.

[KT98] M. Keel and T. Tao, *Endpoint Strichartz estimates*, Amer. J. Math. **120** (1998), no. 5, 955–980.

[Kwo89] Man Kam Kwong, *Uniqueness of positive solutions of $\Delta u - u + u^p = 0$ in $\mathbb{R}^n$*, Arch. Rational Mech. Anal. **105** (1989), no. 3, 243–266.

[Mer93] F. Merle, *Determination of blow-up solutions with minimal mass for nonlinear Schrödinger equations with critical power*, Duke Math. J. **69** (1993), no. 2, 427–454.

[Oh89] Yong-Geun Oh, *Cauchy problem and Ehrenfest's law of nonlinear Schrödinger equations with potentials*, J. Diff. Eq. **81** (1989), no. 2, 255–274.

[Oza91] T. Ozawa, *Long range scattering for nonlinear Schrödinger equations in one space dimension*, Comm. Math. Phys. **139** (1991), 479–493.

[RS75] Michael Reed and Barry Simon, *Methods of modern mathematical physics. II. Fourier analysis, self-adjointness*, Academic Press [Harcourt Brace Jovanovich Publishers], New York, 1975.

[Str74] W. A. Strauss, *Nonlinear scattering theory*, Scattering theory in mathematical physics (J. Lavita and J. P. Marchand, eds.), Reidel, 1974.

[Str77] ———, *Existence of solitary waves in higher dimensions*, Comm. Math. Phys. **55** (1977), no. 2, 149–162.

[Str81] ———, *Nonlinear scattering theory at low energy*, J. Funct. Anal. **41** (1981), 110–133.

[Wei83] Michael I. Weinstein, *Nonlinear Schrödinger equations and sharp interpolation estimates*, Comm. Math. Phys. **87** (1982/83), no. 4, 567–576.

[Yaj96] Kenji Yajima, *Smoothness and non-smoothness of the fundamental solution of time dependent Schrödinger equations*, Comm. Math. Phys. **181** (1996), no. 3, 605–629.

[YZ01] Kenji Yajima and Guoping Zhang, *Smoothing property for Schrödinger equations with potential superquadratic at infinity*, Comm. Math. Phys. **221** (2001), no. 3, 573–590.

MAB, UMR CNRS 5466, Université Bordeaux 1, 351 cours de la Libération, 33 405 Talence cedex, France

*E-mail address*: `carles@math.u-bordeaux.fr`